\documentclass[12pt]{article}
\usepackage{amssymb}
\font\smallit=cmtt10

\newtheorem{theorem}{Theorem}
\newtheorem{lemma}{Lemma}

\begin{document}

\begin{center}
{\bf SUMS OF THE FORM $1/x_1^k + \cdots + 1/x_n^k$ MODULO A PRIME }
\vskip20pt
{\bf Ernie Croot \footnote{Partially Supported by an NSF grant.}} \\
{\smallit Department of Mathematics, Georgia Institute of Technology, 
Atlanta, GA  30332} \\
{\tt ecroot@math.gatech.edu} 
\end{center}
\vskip60pt

\centerline{\bf Abstract}

\noindent Using a sum-product result due to Bourgain, Katz, and Tao,
we show that for every $0 < \epsilon \leq 1$, and every integer $k \geq 1$,
there exists an integer $N = N(\epsilon, k)$, such that for every 
prime $p$ and every residue class $a \pmod{p}$, there exist positive integers
$x_1,...,x_N \leq p^\epsilon$ satisfying 
$$
a\ \equiv\ {1 \over x_1^k} + \cdots + {1 \over x_N^k} \pmod{p}.
$$

\thispagestyle{empty}
\baselineskip=15pt
\vskip30pt

\section*{I.  Introduction}

In the monograph \cite{erdos}, among the many questions asked by Erd\H os
and Graham was the following:  Is it true that for every 
$0 < \epsilon \leq 1$ there exists a number $N$ such that for every
prime number $p$, every residue
class $a \pmod{p}$ can be expressed as 
$a \equiv 1/x_1 + \cdots + 1/x_N \pmod{p}$, where $x_1,...,x_k$ are positive
integers $\leq p^\epsilon$?  This question was answered in the affirmative
by Shparlinski \cite{shparlinski} using a result due to Karatsuba
\cite{karatsuba} (actually, a simplified version of Karatsuba's result, due
to Friedlander and Iwaniec \cite{friedlander}).

A natural question that one can ask, and which Shparlinski recently posed
to me, was whether this result can be extended to reciprocal powers.  
Unfortunately, in this case, the methods of Karatsuba do not give a 
bound on $N$ (at least not using an obvious modification of his
argument).  Fortunately, there is a powerful result due to
Bourgain, Katz, and Tao \cite{bourgain} which can be used to
bound certain exponential sums, and which can be used to solve our 
problem:

\begin{theorem}[Bourgain,Katz,Tao] \label{sum_product} 
Let $A$ be a subset of a finite field
${\mathbb Z}/p{\mathbb Z}$.  If $p^\delta < |A| < p^{1-\delta}$
for some $\delta > 0$, then $|A+A| + |A\cdot A| \geq c|A|^{1+\theta}$, 
where $\theta = \theta(\delta) > 0$ and $c = c(\delta) > 0$. 
\end{theorem}

Using this result, we prove the following theorem, which is just a 
restatement of the problem posed by Shparlinski:

\begin{theorem}  \label{shparlinski_theorem}
For every $0 < \epsilon \leq 1$, and every integer $k \geq 1$, there
exists an integer $N = N(\epsilon,k)$ such that for every prime $p \geq 2$,
and every integer $0 \leq a \leq p-1$, there exist integers $x_1,...,x_N$
such that $1 \leq x_i \leq p^{\epsilon}$, and
$$
a\ \equiv\ {1 \over x_1^k} + \cdots + {1 \over x_N^k} \pmod{p}.
$$
\end{theorem} 

\noindent {\bf Comments.  }  A more general theorem can perhaps be 
proved here, as was suggested to me by Shparlinski in an email.  Basically, 
suppose $S = S(p) \subseteq \{1,...,p-1\}$ is an infinite sequence of sets, indexed by
primes $p$ satisfying the following two conditions
\bigskip

1.  The sets $S(p)$ are multiplicative in the sense that $1 \in S$, and 
if $s,t \in S$ satisfy $st \leq p-1$, then $st \in S$; and,

2.  There exists an absolute constant $0 < \theta \leq 1$ so that for every 
$0 < \epsilon \leq 1$, and $p$ sufficiently large, the set $S(p)$ contains
at least $p^{\epsilon \theta}$ elements $\leq p^\epsilon$.
\bigskip

\noindent Then, there exists an integer
$J = J(\epsilon) \geq 1$ such that for every sufficiently large prime $p$,
and for every residue class $r \pmod{p}$, there exist integers 
$x_1,...,x_J \in S(p)$, all of size at most $p^\epsilon$, such that
$$
r\ \equiv\ {1 \over x_1} + \cdots + {1 \over x_J} \pmod{p},
$$

An example of a set $S = S(p)$ satisfying the properties above, is the 
set of positive integers $\leq p-1$ having no prime divisors greater than $\log^2 p$.
It is well known that the number of elements in this multiplicative set up to
$p^\epsilon$ is at least $p^{\epsilon/2 + o(1)}$ (see, for example, 
\cite{granville})

\vskip30pt

\section*{II.  Proof of Theorem \ref{shparlinski_theorem}}

\vskip30pt

First, we note that it suffices to prove the result only for sufficiently
large primes $p$, as we may enlarge $N = N(\epsilon,k)$ as needed so
that the theorem holds for all prime $p < p_0$, for some $p_0$.  
We also may assume $0 < \epsilon < \epsilon_0(k)$, for any function
$\epsilon_0(k)$ that we might happen to need, since if the conclusion of
the theorem holds for these smaller values of $\epsilon$, then it holds
for any larger value of $\epsilon$.  In fact, we will use 
$\epsilon_0(k) = 1/5k$ in the proof of our theorem. 

Let $0 < \beta < 1/5k$ be some parameter, to be chosen later,
and let $u$ be the largest integer less than $\beta^{-1}/(2k)$, and
consider the set
$$
S\ =\ \left \{ {1 \over p_1^k} + \cdots + {1 \over p_u^k} \pmod{p}\ :\ 
2 \leq p_1 < \cdots < p_u \leq p^\beta,\ p_i\ {\rm prime} \right \},
$$
which will be non-empty for $p$ sufficiently large.  We claim that
\begin{equation} \label{S_size}
|S|\ >\ {p^{1/(2k) - \beta} \over u! \log^u p}
\end{equation}
for $p$ sufficiently large, which would follow from the prime
number theorem if we had that all the sums in $S$ were distinct modulo $p$.
To see that they are, suppose that we had
$$
{1 \over p_1^k} + \cdots + {1 \over p_u^k} \equiv {1 \over q_1^k} + 
\cdots + {1 \over q_u^k} \pmod{p},
$$
where the left and right side of the congruence are elements of $S$,
where the $p_1,...,p_k$ and $q_1,...,q_k$ are increasing sequences.  
Multiplying through
by $(p_1 \cdots p_u q_1 \cdots q_u)^k$ on both sides and moving terms
to one side of the congruence, we get that
\begin{equation} \label{simplified}
\sum_{j=1}^u \left ( (q_1\cdots q_u)^k \prod_{i=1 \atop i \neq j}^u p_i^k 
- (p_1 \cdots p_u)^k \prod_{i=1 \atop i \neq j}^u q_i^k \right ) 
\equiv 0 \pmod{p}.
\end{equation}
Since all the terms in the sum are smaller than 
$p^{(2 u -1)k \beta} < p/u$ (for $p$ sufficiently large), 
we deduce that if (\ref{simplified}) holds, then 
$$
\sum_{j=1}^u \left ( (q_1\cdots q_u)^k \prod_{i=1 \atop i \neq j}^u p_i^k 
- (p_1 \cdots p_u)^k \prod_{i=1 \atop i \neq j}^u q_i^k \right ) 
\ =\ 0;
$$
and so,
$$
{1 \over p_1^k} + \cdots + {1 \over p_u^k}\ =\ {1 \over q_1^k} + 
\cdots + {1 \over q_u^k}.
$$ 
It is obvious then that the $p_i = q_i$, and (\ref{S_size}) now follows.
\bigskip

Let $S_0 = S$, and consider the sequence of subsets of ${\mathbb Z}/p{\mathbb Z}$,
which we denote by
$S_1,S_2,...$, where 
$$
S_{i+1}\ =\ \left \{ \begin{array}{rl} S_i + S_i,\ &{\rm if\ }|S_i+S_i| > |S_i S_i|;\ {\rm and} \\
S_iS_i,\ &{\rm if\ }|S_i+S_i| \leq |S_iS_i|. \end{array} \right .
$$
We continue constructing this sequence until we reach the set $S_n$ satisfying
\begin{equation} \label{Sn_target}
|S_n| > p^{2/3}.
\end{equation}
Using Theorem \ref{sum_product} we can produce
a non-trivial upper bound on the size of $n$ for $\beta < 1/5k$:  
Let $\delta = 1/4k$,
and let $c = c(\delta)$, $\theta = \theta(\delta)$ be as in Theorem
\ref{sum_product}.  Then, for $p$ sufficiently large, we will have
$$
p^\delta\ <\ |S_0|\ =\ |S|\ <\ p^{1-\delta},
$$
and the same inequality will hold for $S_1,S_2,...,S_{n-1}$.  Now, applying 
Theorem \ref{sum_product}, we deduce that
$$
|S_{i+1}|\ >\ c |S_i|^{1 + \theta};
$$
and so, 
$$
|S_1|\ >\ c |S_0|^{1+\theta},\ \ {\rm and\ for\ }j=2,...,n,\ 
|S_j|\ >\ c^{1 + (1+\theta)^{j-1}} |S_0|^{(1+\theta)^j}.
$$
From this inequality and (\ref{Sn_target}) , we deduce that
$$
 \left ( {1 \over 2k} - \beta \right ) (1+\theta)^n + o(1)\ >\ {2 \over 3},
$$
where the $o(1)$ tends to $0$ as $p$ tends to infinity; and
so, since $\beta < 1/5k$, our sequence $S_0,S_1,...,S_n$ finishes with
$$
n\ <\ {\log(3k) \over \log(1+\theta)} + 1
$$
for $p$ sufficiently large.

Now, every element of $S_0$ is a sum of at most $u$ terms;
each element of $S_1$ is a sum of at most $u^2$ terms; and, by an induction
argument, each element of $S_n$ is a sum of at most $u^{2^n}$ terms.  
Also, each element of $S_n$ is a sum of terms of the form 
$1/q_1^k \cdots q_{2^n}^k$, where $q_1\cdots q_{2^n} \leq p^{2^n\beta}$.

Now, let $\beta = \epsilon/2^{n+1}$.  If $\epsilon < 1/5k$, then this value of 
$\beta < 1/5k$ (recall we said that $\epsilon$ is allowed to be bounded from
above by a function of $k$).   Let $h$ 
$$
h\ =\ u^{2^{\left \lceil {\log(3k) \over \log(1 + \theta)} \right \rceil}},
$$ 
and define 
$$
T\ =\ \left \{ {1 \over q_1^k} + \cdots + {1 \over q_h^k}\ \pmod{p}\ :\ 
2 \leq q_1,...,q_h \leq p^{\epsilon/2} \right \}
$$
Here, $q_1,...,q_h$ are not restricted to being prime numbers.
Since $|T| \geq |S_n|$, we have that $|T| > p^{2/3}$ for $p$ sufficiently large.

Now we use the following simple lemma, which has appeared in many works
before, and uses a standard bilinear exponential sums technique:

\begin{lemma}  Suppose that $T \subseteq {\mathbb Z}/p{\mathbb Z}$
satisfies $|T| > p^{1/2 + \beta}$.  Then, every residue class modulo $p$ 
contains an integer of the form $x_1+\cdots + x_J$, where the 
$x_1,...,x_J$ are all of the form $t_1t_2$, where $t_1,t_2 \in T$, and where
$J = \lfloor 2(1+2\beta)/\beta \rfloor  + 1$.
\end{lemma}

\noindent {\bf Proof of the Lemma.}  First, we consider the exponential sums
$$
h(a)\ =\ \sum_{t \in T} e \left ( {at \over p} \right ),
$$
and 
$$
f(a)\ =\ \sum_{t_1,t_2 \in T} e \left ( {at_1t_2 \over p} \right ).
$$
We have from Parseval's identity and the Cauchy-Schwarz inequality that 
for $a \not \equiv 0 \pmod{p}$,
\begin{eqnarray}
|f(a)|\ &\leq&\ \sum_{t_1 \in T} \left | \sum_{t_2 \in T} e \left ( {a t_1 t_2 \over p} \right )
\right | \nonumber \\
&\leq& \left ( \sum_{t_1 \in T} 1 \right )^{1/2} \left ( \sum_{t_2 \in T} 
|h(at_2)|^2 \right )^{1/2} \nonumber \\
&=&\ p^{1/2} |T|\ \leq\ |f(0)|^{(1+\beta)/(1+2\beta)}. \nonumber
\end{eqnarray}  
Now, if we let $J$ be the least integer greater than 
$$
2\left (1 - {1+\beta \over 1+2\beta} \right )^{-1}\ =\ 
{2 (1+2\beta) \over \beta},
$$
then we have that for $a \not \equiv 0 \pmod{p}$, 
\begin{eqnarray}
|f(a)|^J\ &<&\ |f(0)|^{ J (1+\beta)/(1+2\beta)}\ <\ |f(0)|^J |f(0)|^{-J \beta/(1+2\beta)}
\nonumber \\
&\leq&\ |f(0)|^{J-2}\ \leq\ {|f(0)|^J \over p}. \nonumber 
\end{eqnarray}
Thus, given an integer $r$, the number
\begin{eqnarray}
&& \#(x_1,...,x_J\ :\ x_i = t_1t_2;\ t_1,t_2 \in T;\ {\rm and\ } x_1 + \cdots + x_J 
\equiv r \pmod{p} )\nonumber \\
&&\ \ \ \ \ =\ {1 \over p} \sum_{a=0}^{p-1} f(a)^J e(-ar/p) 
\ \geq\ {f(0) \over p}\ -\ {1 \over p} \sum_{1 \leq a \leq p-1} |f(a)|^J \nonumber \\
&&\ \ \ \ \ \geq\ {f(0) \over p}\ -\ {(p-1) f(0) \over p^2}
\ >\ 0. \nonumber
\end{eqnarray}
This proves the lemma.\ \ \ \ \ \ \ $\blacksquare$
\bigskip

From this lemma, we deduce that for every residue class $r$ modulo $p$,
there exist integers $t_1,...,t_{16}$, such that
$$
r\ \equiv\ t_1 t_2 + t_3t_4 + \cdots + t_{15} t_{16} \pmod{p},
$$
where $t_1,...,t_{16} \in T$.  This sum can be expressed as a sum of at 
most $16h^2$ terms of the form $1/(q q')^k$, where $q,q' < p^{\epsilon/2}$.
This then proves the theorem, since $h$ depends only on $k$ and
$\epsilon$.

\section*{\normalsize Acknowledgements}

I would like to thank Igor Shparlinski for the interesting problem, which is
the one solved by Theorem \ref{shparlinski_theorem}.


\begin{thebibliography}{999}

\bibitem{bourgain} J. Bourgain, N. Katz, and T. Tao, 
{\it A Sum-Product Estimate in Finite Fields, and Applications},
Preprint on the Arxives.

\bibitem{erdos} P. Erd\H os and R. L. Graham, {\it Old and New Problems and
Results in Combinatorial Number Theory}, Univ. Gen\` eve, Geneva, 1980.

\bibitem{friedlander} J. Friedlander and H. Iwaniec, {\it Analytic Number Theory 
(Kyoto, 1996)}, Cambridge University Press, Cambridge, 1997.

\bibitem{granville} A. Granville, {\it Smooth Numbers:  Computational Number Theory
and Beyond}, MSRI Workshop Notes.

\bibitem{karatsuba} A. A. Karatsuba, {\it Fractional Parts of Functions of
a Special Form}, Izv. Ross. Akad. Nauk Ser. Mat. {\bf 59} (1995), 61-80. 

\bibitem{shparlinski} I. Shparlinski, {\it On a Question of Erd\H os and
Graham}, Arch. Math. (Basel) {\bf 78} (2002), 445-448.

\end{thebibliography}
\end{document}